\documentclass[sigplan,screen]{acmart}

\setcopyright{cc}
\setcctype{by}
\acmDOI{10.1145/3703595.3705885}
\acmYear{2025}
\copyrightyear{2025}
\acmISBN{979-8-4007-1347-7/25/01}
\acmConference[CPP '25]{Proceedings of the 14th ACM SIGPLAN International Conference on Certified Programs and Proofs}{January 20--21, 2025}{Denver, CO, USA}
\acmBooktitle{Proceedings of the 14th ACM SIGPLAN International Conference on Certified Programs and Proofs (CPP '25), January 20--21, 2025, Denver, CO, USA}
\received{2024-09-08}
\received[accepted]{2024-11-19}

\usepackage[utf8]{inputenc}
\usepackage[T1]{fontenc}
\usepackage{xspace}
\usepackage{cleveref}
\usepackage{microtype}

\usepackage[vcentermath]{genyoungtabtikz}

\usepackage{listings}
\definecolor{dkblue}{rgb}{0,0.1,0.5}
\definecolor{lightblue}{rgb}{0,0.5,0.5}
\definecolor{dkgreen}{rgb}{0,0.4,0}
\definecolor{dk2green}{rgb}{0.4,0,0}
\definecolor{dkviolet}{rgb}{0.6,0,0.8}
\definecolor{brick}{rgb}{0.6,0.2,0.25}

\let\verb=\lstinline

\usepackage{commath}
\newcommand{\Sage}{\texttt{SageMath}\xspace}
\newcommand{\Coq}{\texttt{Coq/Rocq}\xspace}

\newcommand{\MC}{\texttt{MathComp}\xspace}
\newcommand{\LR}{Littlewood-Richardson\ }
\newcommand{\RS}{Robinson-Schensted\ }
\newcommand{\var}[1]{\text{\verb{#1}}}

\newcommand\ie{\textit{i.e.}\xspace}
\newcommand\eg{\textit{e.g.}\xspace}

\newcommand{\N}{{\mathbb N}}

\newcommand{\K}{{\mathbb K}}
\newcommand{\SG}{{\mathfrak S}}
\newcommand{\std}{\operatorname{Std}}
\newcommand{\sgn}{\operatorname{sgn}}

\newcommand{\sym}{\mathrm{sym}}

\newcommand{\partof}{\vdash}                    
\newcommand{\shape}{\operatorname{shape}} 

\newcommand{\red}[1]{{\color{red} #1}}
\newcommand{\grn}[1]{{\color{green} #1}}
\newcommand{\blu}[1]{{\color{blue} #1}}

\newcommand{\alphX}{{\mathbb X}}

\newcommand{\alphA}{{\mathbb A}}
\newcommand{\placteq}{\equiv_\text{Plact}}

\newtheorem{problem}{Problem}
\newtheorem{DEFN}{Definition}

\newtheorem{ALGO}{Algorithm}

\begin{abstract}
  We present a library of formalized results around symmetric functions and
  the character theory of symmetric groups. Written in Coq/Rocq and based on
  the Mathematical Components library, it covers a large part of the contents
  of a graduate level textbook in the field. The flagship result is a proof of
  the Littlewood-Richardson rule, which computes the structure constants of
  the algebra of symmetric function in the schur basis which are integer
  numbers appearing in various fields of mathematics, and which has a long
  history of wrong proofs. A specific feature of algebraic combinatorics is the
  constant interplay between algorithms and algebraic constructions:
  algorithms are not only in computations, but also are key ingredients in
  definitions and proofs. As such, the proof of the Littlewood-Richardson rule
  deeply relies on the understanding of the execution of the Robinson-Schensted
  algorithm. Many results in this library are effective and actually used in
  computer algebra systems, and we discuss their certified implementation.
\end{abstract}


\begin{document}

\title{Machine Checked Proofs and Programs in Algebraic Combinatorics}
\author{Florent Hivert}
\email{florent.hivert@lisn.fr}
\orcid{0000-0002-7531-5985}
\affiliation{%
  \institution{University Paris-Saclay}
  \city{Orsay}
  \country{France}
}
\affiliation{%
  \institution{LISN, LMF, CNRS, INRIA}
  \city{Orsay}
  \country{France}
}

\maketitle

\section{Introduction : algebraic combinatorics}

Algebraic combinatorics is an area of mathematics that studies the
interactions of algebra (\eg linear algebra, polynomials, symmetry groups and
representation theory) with algorithms (data structures, searching and
sorting). These interactions go both ways. One can use abstract
algebra, in particular group theory and representation theory to study
combinatorial objects (\eg isomorphism of combinatorial structure
such as graph). Conversely, many combinatorial techniques can apply to
algebra. Algorithms are used not only to compute but to construct algebraic
interpretations and to prove identities.

In order to illustrate the interplay between data-structures and algebra, let us start with a
simple yet deep example: the coefficients of the expansion of a fully factored
polynomial are given by
\begin{equation}
  \prod _{{j=1}}^{n}(X -\lambda_{j})=X
  ^{n}-e_{1}X ^{{n-1}}+
  e_{2}X ^{{n-2}}+
  \cdots +(-1)^{n}e_{n},
\end{equation}
where the $e_k$ are called the \emph{elementary symmetric polynomials} and are defined by formula
\begin{equation}
  e_k := \sum_{\abs{S} = k} \prod_{i\in S} \lambda_i
\end{equation}
where the sum ranges over subsets $S$ of $\{1,\dots,n\}$ of cardinality
$k$. For each $k$, the number of such subsets is the binomial coefficient
$\binom{n}{k}$ appearing in Newton's binomial formula
\begin{equation}
(a+b)^n = \sum_{k=0}^{n} \binom{n}{k} a^k b^{n-k}\,.
\end{equation}
In these two toy examples, the key of the \emph{algebraic} formula is the
\emph{combinatorial} enumeration of the subsets of
$\{1,\dots,n\}$. Interestingly enough, elementary symmetric functions are a
special case of Schur functions, themselves the crux of the \LR rule.

\smallskip Computer experimentation largely fuels research in algebraic
combinatorics, a quite distinctive feature of this area of mathematics. The
community has heavily invested in the computer algebra systems such as
Sagemath~\cite{Sagemath}. To give a figure, the size of the \verb|combinat|
directory of the \Sage software comes close to half a million lines of code
(including tests and documentation). However, our involvement in that project
showed this huge amount of code exhibits many bugs, lots of them being
specifically on the particular cases (empty domains, zero case, etc.). The
common answer to this issue is to deploy tests. And indeed, tests are very
relevant in combinatorics, as it is often possible to expand the full list of
possible inputs up to some size. In our experience this works relatively well.

Unfortunately, a second specific feature of the field is that algorithms are
at the core of \emph{reasoning}, as, for example, in the case of bijective
proofs: two finite sets $A$ and $B$ are proved to be
equinumerous by introducing two functions $f:A\to B$ and $g:B\to A$, and by
proving that they are inverse of each other:
\[
  \forall a\in A, g(f(a)) = a
  \quad\text{and}\quad
  \forall b\in B, f(g(b)) = b.
\]
Very often, $f$ and $g$ are defined by complicated algorithms and quite often
their implementation is wrong, which sheds some doubt about the correctness of
proofs.

\medskip
The Littlewood-Richardson is a famous case of problematic combinatorial proof.
In 1934, D. E. Littlewood and A. R. Richardson described a rule to compute
some very important non-negative integer constants, by counting the number of
solutions to an elementary, yet intricate combinatorial problem. The first
correct proof had to wait for the mid 70's and several wrong proofs were
published and widely accepted by the community:
\begin{quotation}[Wikipedia]
  The Littlewood–Richardson rule is notorious for the number of errors that
  appeared prior to its complete, published proof. Several published attempts
  to prove it are incomplete, and it is particularly difficult to avoid errors
  when doing hand calculations with it: even the original example in
  D. E. Littlewood and A. R. Richardson (1934) contains an error.
\end{quotation}
In his review of James and Kerber book~\cite{JamesKerber.1984}, Jacob Towber
wrote~\cite{Towber83}:
\begin{quotation}
  It seems that for a long time the entire body of experts in the field was
  convinced by these proofs; [...]  How was it possible for an incorrect proof
  of such a central result in the theory of the symmetric groups to have been
  accepted for close to forty years? The level of rigor customary among
  mathematicians when a combinatorial argument is required, is (probably quite
  rightly) of the nonpedantic hand-waving kind; perhaps one lesson to be drawn
  is that a higher degree of care will be needed in dealing with such
  combinatorial complexities as occur in the present level of development of
  Young's approach.
\end{quotation}
Our initial motivation for the present work was thus to evaluate whether a fully
computer-checked formalization of this result was feasible at all, and if so
to evaluate its cost.

\medskip
We thus formalized a consistent corpus of basic results in the
field of symmetric functions and characters of the symmetric group. The
formalization covers most of the results of Sagan's book ``The Symmetric
Group: Representations, Combinatorial Algorithms, and Symmetric
 Functions''~\cite{Sagan01},
in the collection ``Graduate Texts in Mathematics'',
the algebraic part being equivalent to the widely cited first chapter of
Macdonald's much more advanced book~\cite{Macdonald.1995}. The full
formalization extends along 35k lines of Coq/Rocq code, among which only 5k to
10k are specialized results aimed at the proof of the \LR rule, depending on
whether we count the study of the \RS correspondence and the plactic monoid
among those. In particular, most of the basic results regarding various
combinatorial objects, the combinatorics of the symmetric group and the
elementary theory of symmetric functions are largely reusable. Some of the
side results presented in \cref{par.side-results} are also of general
interest.

The \LR rule is the flagship result of the library. We chose to formalize the
original proof, due to Schützenberger~\cite{SchutzLR}, rather than a modern shorter
proof. The first reason is that some eminent colleague expressed doubts about the lack of gap in this proof. But most importantly, the proof
involves a very wide and representative spectrum of techniques used in the
fields. The proof actually features intricate and quite subtle combinatorial games involving basic enumerations, partial ordering, algorithms and invariants,
rewriting, advanced bijections, linear algebra, or group manipulations.

To the best of our knowledge, this is the first time such a large corpus has
been formalized. All the formal development about symmetric polynomials we are
aware of such as~\cite{Symmetric_Polynomials-AFP,Cohen-Djalal,LeanSymmetric}
are mostly limited to two very basic results: the fundamental theorem of
symmetric polynomial stating that any symmetric polynomial can be expressed a
polynomial in the elementary ones and Newton's identities expressing
recursively the elementary one on power sums. The proof of the latter identity
takes around 200 lines among the 2500 of the file \texttt{sympoly.v} devoted
to the classical bases of symmetric polynomials.
\smallskip

The proof was formalized in \Coq over the Mathematical Components library
(\MC). The main reason why we choose \MC as a basis was because we were
convinced by Gonthier's insights about ``Mathematical engineering'' and the
``programs-as-proof'' style of designing libraries~\cite{Gonthier4col,
  MathCompBook}. The small scale reflection methodology, where most
definitions are written using computational boolean predicates, proved very
well-suited to combinatorics, as it saves the user the pain of dealing with
myriads of particular cases, most of that are trivial. Remarkably, those are
the very same cases where many bugs are found in implementations.  We feel
that having the system deal with this bureaucracy automatically is
instrumental for formalizing combinatorics. Finally, being able to \emph{test}
executable combinatorial definitions, that are very easy to get wrong, was
also quite useful.  Small scale reflection has proven very effective to reduce

Another pleasant outcome of this methodology has been  that even our first written formal statement of the \LR rule could compute a few small numbers. We nevertheless provide a more efficient
backtracking algorithm. The later is close to those implemented in actual computer algebra
systems, such as Ander Buch's library~\cite{lrcalc}, albeit being a few orders of magnitude
slower (mostly due to the lack of mutable data structures). Note that there is some very strong
complexity theoretical evidence that there is no fundamentally better way to
compute these numbers (see \cref{subsec:comput}).

A contribution of this library is thus the first certified implementation of
the \LR rule.  We anyway emphasize that the main contribution of this work is
not about proving the correctness of some algorithms, but rather about
formalizing mathematical statements asserting that some algebraic constants
can be computed by counting the number of solutions of a combinatorial game.

\subsection*{Outline of the paper}

The paper starts by introducing the mathematical context of symmetric Schur
functions and the \LR coefficients (\cref{sec:context}) together with some
comments concerning their appearance in various areas, as well as the
relevance of the rule as a practical mean to compute those. After a short
section introducing the terminology and related formalization of very basic
language theory (\cref{sec:language}), we describe the main combinatorial
objects, namely partitions, tableaux and Yamanouchi words
(\cref{sec:part-tab}).

We then switch to the algebraic background (\cref{sec:sf}) of symmetric
functions culminating with the statement of the \LR rule
in~\cref{theo:LR-rule} and its formalized statement. The next section
(\cref{sec:SymGroup}) is a brief description of the translation of the result
on symmetric function in the language of character of the symmetric groups
through Frobenius characteristic isomorphism. The next section is a broad
description of the content of the library~(\cref{sec:overview}).

The core of the paper is in \cref{sec:LRproof}, where we give some insight
about the actual proof of the LR rule. We also explain how
\cref{algo:Schensted}, due to Schensted, computes some invariants and normal
forms of a quotient of the free monoid by some simple rewriting rule called
the plactic monoid~(\cref{theo:plact}). We then briefly describe the central
idea, which is to lift the computation at a noncommutative level
(\cref{subsec:noncomm}). We further discuss the implementation
(\cref{sec:implem}) and conclude with some comments on our development. The
paper is extended by an appendix which is the table of contents of the
library.

\section{Mathematical context}
\label{sec:context}
\subsection{Symmetric functions}
The ring of \emph{symmetric functions}~\cite{Macdonald.1995} is defined as the
(categorical) limit when $n$ goes to infinity of the ring of symmetric
polynomials in $n$ indeterminates. This ring allows to express relations
between symmetric polynomials independently from the number of indeterminates
involved, although its elements are neither polynomials nor functions. In
particular, this ring plays an important role in representation theory of the
symmetric group and of the general linear group, and in various geometry
problems (see next section for some examples).

A particular linear basis $(s_\lambda)_\lambda$ of this algebra, whose elements are called the \emph{Schur functions}, plays a key role in these various contexts. Recall that linear bases of symmetric functions are indexed by
\emph{integer partitions} $\lambda$, i.e., by non-increasing sequences of positive
integers. Actually, the ring of symmetric functions can be endowed with a structure of algebra, and the product of two Schur functions
can itself be expressed as a linear combination of Schur functions (see later for
examples):
\begin{equation}
  s_\lambda s_\mu = \sum_{\nu} C_{\lambda, \mu}^{\nu}\ s_\nu\,.
\end{equation}
Coefficients $C_{\lambda, \mu}^{\nu}$ in the decomposition are nonnegative integers called 
\emph{\LR coefficients}. The \emph{\LR rule} describes these coefficients as the number of certain combinatorial configurations called \emph{\LR tableaux}.
The precise definition of these configurations is rather intricate and
intertwines several types of constraints on the way to fill a finite
collection of boxes, arranged according a prescribed shape, with numbers. It
makes the rule difficult to state, to use and even more to prove. As already
mentioned, the first correct proof, by Schützenberger~\cite{SchutzLR}, was
given in 1977, four decades after the rule was first stated. This proof
involves numerous combinatorial ingredients, and we already mentioned that to
some combinatorialists considered this proof as ``somewhat gappy''. The
present work roughly follows this original proof, as presented
in~\cite{Lothaire.2002-Plact}. The algebraic step is however closer to a later
alternate presentation~\cite{NCSF6}. It shows in particular that this original
proof actually has no crucial gap. Note that Schützenberger's orginal proof
was followed by many attempts at simplifying the
argument~\cite{Zelevinsky81,NCSF6,VanLeeuwen01,Stembridge02}. The
combinatorial study of \LR coefficients is by the way still a very active
research topic~(see for example~\cite{KnutsonTao99, Ikenmeyer2017, Rosas2023}
to only name a few).

Beyond validating this original proof, our motivation for a formal study of the \LR rule is actually two-fold. First, \LR coefficients are ubiquitous in mathematics but also physics and
chemistry. Second, there is in some sense no better way to compute them than
combinatorial rules of the same. The rest of this section elaborates on these two points.

\subsection{Various views on \LR coefficients}
\begin{quotation}[Gordon James (1987)]
  Unfortunately the Littlewood–Richardson rule is much harder to prove than
  was at first suspected. The author was once told that the
  LR rule helped to get men on the moon but was not proved
  until after they got there. The first part of this story might be an
  exaggeration.
\end{quotation}
Although most certainly a joke, this quotation from an expert in representation theory is meant to stress how ubiquitous those
numbers are. We mention here a few appearances only, and refer the interested reader
 to Fulton's overview~\cite{Fulton96} or Macdonald's
book~\cite{Macdonald.1995}.
\begin{itemize}
\item Coefficients $C_{\lambda\mu}^\nu$ can equivalently be defined as the inner product
  $\langle s_\nu \mid s_\lambda s_\mu\rangle$. They also provide the expansion
  of a Schur function on the union of two disjoint sets of variables; 
\item They count the multiplicity of induction or restriction of irreducible
  representations of the symmetric group;
\item By Schur-Weyl duality, they also count the multiplicity of the tensor
  product of the irreducible representations of linear groups or special
  linear groups;
\item They also have various geometrical interpretations: for example, they
  are the intersection number in a grassmanian variety and also appear in
  the cup product of its cohomology;
\item They are connected to the Horn problem about the eigenvalues of the sum
  of two matrices with prescribed eigenvalues.
\item They are related to extensions of abelian $p$-groups (\ie, whose order
  is a power of a prime number $p$) through the Hall algebra.
\item Finally, their connections with finite group theory explain their role in
  quantum physics (orbitals and spectrum of
  atoms~\cite{Wybourne1992}) and in chemistry~\cite{Chemistry}.
\end{itemize}

\subsection{Complexity}
\label{subsec:comput}

From an effective point of view, such a rule is the only good way to
compute these numbers, as the
computation of the \LR coefficients is
$\#P$-complete~\cite{Narayanan06}. Recall that $\#P$ is the complexity class
of counting problem (\ie with an answer in $\N$) analogue to the complexity class
$NP$ for decision problem (\ie with a boolean answer). This roughly means
that we cannot hope for a significantly better algorithm for computing these numbers 
than enumerating the solutions of a combinatorial problem such as
\LR tableaux.

Note that other combinatorial models exist that calculate these coefficients,
such as Knutson and Tao honeycombs~\cite{KnutsonTao99}. Note as well that \LR
coefficients and honeycombs appear in Mulmuley's geometric complexity
theory~\cite{Mulmuley01,Mulmuley12}, a research program to prove that
$P\neq NP$ based on invariant theory and algebraic geometry.


\section{Alphabet, words and languages}
\label{sec:language}
In this paper, an \emph{alphabet} $\alphA$ is a totally ordered set. It is
thus modeled as an instance of the \MC \verb|orderType| structure.  For
technical reasons, we assume that $\alphA$ is inhabited (non empty), \ie an
instance of \verb|inhOrderType|. A \emph{word} is a finite sequence of
\emph{letters} (elements of $\alphA$), i.e., an inhabitant of type
\verb|seq Alph| in \MC. The empty word \verb|[]| and the concatenation
\verb|u ++ v| of words \verb|u| and \verb|v| endow the set of word $\alphA^*$
with a structure of monoid, which is actually the free monoid. We use the
general convention that symbols $a,b,c$ denote letters in $\alphA$ and $u,v,w$
denote words, in $\alphA^*$. The main concrete alphabets that are used in the
library are the integers, formalized by type \verb|nat|, and the the finite
set $\{0,\dots, n-1\}$, formalized by the type \verb|'I_n| of finite ordinals.

As a side note, the first versions of the project contained an extensive
library for defining and relating structures of (partially or totally, inhabited) ordered types and sequences over an ordered type.  The hierarchy involved many diamonds,
which made it very hard to correctly set up inheritance relations. We actually never got it perfectly
correct. It turns out that the Hierarchy Builder~\cite{Cohen20} tool greatly simplified this part of the infrastructure. In particular, the treatment of an inhabited variant of the general case now only requires a hundred of lines of code.

A \emph{language} is just a set of words over an alphabet. At the time when
the library was written, \MC dealt exclusively with finite sets (type
\verb|{set T}|) with a finite underlying type \verb|T| (a structure called
\verb+finType)+~\cite[Chapter~7]{MathCompBook}. Fortunately, for our purpose
it is sufficient to only deal with \emph{homogeneous} languages, that is,
languages of words that all have the same length. Moreover, if the type
\verb|Alph| (of letters) is finite, \MC provides a canonical instance of
\verb|finType| structure on type \verb|d.-tuple Alph|, that is the type of
sequences of elements of \verb|Alph| with prescribed length \verb|d|. We
therefore naturally represent $d$-homogeneous languages as values of type
\verb|{set d.-tuple Alph}|.

\section{Partitions and tableaux}
\label{sec:part-tab}
\subsection{Partitions}

Bases of the algebra of symmetric functions are indexed by (integer) partitions.
\emph{Partitions} are defined as the different ways of decomposing an integer
$n\in\N$ as a sum. For example, there are $7$ partition of $5$, namely
\[ 5=4+1=3+2=3+1+1=2+2+1=2+1+1+1=1+1+1+1+1\,. \] Two decompositions that
differ only by their order are considered equal. To ensure unicity of the
machine representation, we follow the usual convention to sort the summand
in decreasing order.
\begin{DEFN}
  A \emph{partition} $\lambda$ of an integer $n$ is a finite decreasing
  sequence of positive integers
  $(\lambda_0\geq\lambda_1\geq\dots\geq\lambda_{l-1} > 0)$ whose sum is
  $n$. We denote $|\lambda| := n = \lambda_0+\lambda_1+\dots+\lambda_{l-1}$
  the sum and $\ell(\lambda) := l$ the length of the partition $\lambda$. We also denote $\lambda\partof
  n$ the fact that $\lambda$ is a partition of $n$.

  Conventionally, there is only one partition of the integer $0$, namely
  the empty sequence, so that $\lambda\partof0$ means that $\lambda = ()$.
\end{DEFN}
We represent partitions by terms of type \verb+seq nat+.
Equality of partition therefore coincides with Leibnitz equality of
the corresponding sorted lists.  Following the \MC
methodology~\cite[Chapter~5]{MathCompBook}, we introduce a computation 
definition for the predicate \verb|is_part : seq nat -> bool|, using a recursive function.
We also provide various alternative statements in sort \verb|Prop|, whose equivalence 
to the computational variant is stated using reflection lemmas, like \verb+is_partP+:
\begin{lstlisting}
  Fixpoint is_part sh :=     (* Boolean Predicate *)
    if sh is sh0 :: sh'      (* SSR syntax for pattern matching *)
    then (sh0 >= head 1 sh') && (is_part sh')
    else true.
  Lemma is_partP sh : reflect (* Boolean reflection lemma *)
    (last 1 sh != 0 /\ forall i, (nth 0 sh i) >= (nth 0 sh i.+1))
    (is_part sh).
\end{lstlisting}
where \verb|last v s| and \verb|nth v s i| are respectively the last and
$n$-th entry of \verb|s| or \verb|v| if the list is too short or even empty.

Then, we model the finite set $P_n$ of partitions of $n$ by defining a
dependent type \verb|intpartn n|, pairing a sequence and a proof that 
it is a partition of \verb+n+. We provide this type with a canonical instance of 
finite type (\verb+fintype+) structure.  Note that by design choice, most of
the lemmas on partitions (and other combinatorial object such as tableaux or
Yamanouchi words) usually quantify over terms sequences of integers, and feature
a separate assumption that this sequence is a partition, rather than
on dependent pairs of a data and a proof. Dependent types, such as \verb|intpartn n|, 
are actually mostly used to manifest the finiteness of a collection, chiefly for expressing properties about cardinalities.

It is customary to depict a partition by a diagram of boxes called its Ferrers
diagram. Namely the Ferrers diagram of a partition $\lambda := (\lambda_0,
\lambda_1,\dots,\lambda_{l-1})$ is obtained by piling left justified rows of
boxes of respective length $\lambda_0, \lambda_1,\dots,\lambda_{l-1}$. We use
the French convention that puts the longest row at the bottom of the
picture (English literature usually draws them upside down). For example,
\[(7,5,3,2,2)\quad\text{ is depicted as }\quad \yngs(0.5, 2,2,3,5,7).\]

Partitions are partially ordered by the inclusion of their diagrams,
which is formalized by:
\begin{lstlisting}
Definition included inn out :=
  size inn <= size out /\ forall i, nth 0 inn i <= nth 0 out i
\end{lstlisting}
A skew partition is the difference of two included partitions:
\[(7,5,3,2,2) / (4,2,1)\quad\leftrightarrow\quad \gyoungs(0.5,\ \ ,\ \ ,:;\ \
,::;\ \ \ ,::::;\ \ \ ).\] We did not define a specific type for
skew partitions in \Coq: when a skew partition is required we just pass
proofs of \verb+is_part inn+, \verb+is_part out+ and \verb+included inn out+.

\subsection{Tableaux}

Tableaux were invented by Alfred Young to understand the representation of the
symmetric groups.
\begin{DEFN}
  A \emph{Young tableau} or \emph{tableau} for short on $\alphA$ is a filling
  $T$ with letters from $\alphA$ of the diagram of a partition $\lambda$, that
  is, \emph{non decreasing along rows} and \emph{strictly increasing along
    columns}. The partition $\lambda$ is called the \emph{shape} of $T$.

  A \emph{standard tableau} is a tableau over the integer such that each
  integer between $0$ and $n-1$, where $n$ is the sum of the parts of the shape,
  appear only once.

  A \emph{skew tableau} is a tableau whose shape is a skew shape.
\end{DEFN}
Note that in the literature, standard tableau are usually labeled from $1$ to
$n$ instead of $0$ to $n-1$. We made this choice as it was consistent with
\MC type for ordinal (\verb|'I_n|) and permutations. Here is an example of a
tableau, a standard tableau and a skew tableau.
\[
  \Yboxdim{10pt}\scriptsize
  \young(ff,cdd,bccdf,aabeefgh)\qquad
  \young(7,4,258,01369)\qquad
  \gyoung(12,:;00,:::;1,:::;00)\qquad
\]
Following~\cite{Lothaire.2002-Plact}, we formalize tableaux by defining a
\emph{dominance} relation between two consecutive nondecreasing sequence
called \emph{rows}:
\begin{DEFN}
  A non decreasing word $v \in \alphA^*$ is called a \emph{row}. Let $u = x_0
  \dots x_{r-1}$ and $v = y_0 \dots y_{s-1}$ be two rows. We
  say that \emph{$u$ dominates $v$} if $r\leq s$ and for $i = 0,\dots,r-1$,
  $x_i > y_i$.

  A \emph{tableau} is a sequence of non empty rows that is reverse sorted for
  the dominance (strict) order.
\end{DEFN}
As for partitions, we do not introduce a specific type for rows and tableaux but 
simply use type \verb{seq (seq T)}. Tableaux are characterized by
a boolean predicate \verb{is_tableau}, along with its equivalent in
\verb|Prop|. The definition of \verb|dominate u v| is:
\begin{lstlisting}
size u <= size v /\ forall i, i < size u -> nth Z u i > nth Z v i
\end{lstlisting} 
and the one of \verb|is_tableau u v| is:
\begin{lstlisting}
forall i, i < size t -> (nth [::] t i) != [::] (* no empty rows *)
/\ forall i, sorted (@leq_op T) (nth [::] t i)
/\ forall i j, i < j -> dominate (nth [::] t j) (nth [::] t i)
\end{lstlisting}
The \emph{row reading} of a tableau is a crucial notion. It is just the word
obtained from the natural reading (top to bottom and left to right) of a
tableau. For example, the reading of the first tableau above is:
$ffcddbccdfaabeefgh$. It is formally defined by the reversing and then flattening a tableau \verb|t|: \verb|flatten (rev t)|.

Many proofs involve surgery on tableaux. For example, we define a function
\verb{join_tab} that glues the lines of two (possibly skew) tableaux, as
exemplified below:
\begin{lstlisting}
Definition join_tab s t :=
  [seq r.1 ++ r.2 | r <- zip (pad [::] (size t) s) t].
\end{lstlisting}
\[
\var{join_tab}\left(
  \Yboxdim{10pt}\scriptsize
  \ \young(,c,bcc,aab)\,,\ 
  \gyoung(ff,:;dd,:::;df,:::;eefg)\ \right)\ =\
  \Yboxdim{10pt}\scriptsize
  \young(ff,cdd,bccdf,aabeefg)
\]
As an example of surgery, the following lemma asserts that if all the entries
of a tableau \verb{s} are smaller than all the entries of a skew tableau
\verb{t} whose inner shape is the shape of \verb{s}, then the join of
\verb{s} and \verb{t} is itself a tableau:
\begin{lstlisting}
Lemma join_tab_skew s t :
  all (allLtn (to_word s)) (to_word t) ->
  is_tableau s -> is_skew_tableau (shape s) t ->
  is_tableau (join_tab s t).
\end{lstlisting}
Interestingly, though any combinatorialist would consider that such a lemma is
so trivial that is doesn't require any proof, its formal proof requires a 58~line script.

\subsection{Yamanouchi words}

The data structure of Yamanouchi words is equivalent to that of standard tableaux (meaning that they are related by a simple bijection). Yet since the former are one-dimensional, they are easier to manipulate in some cases. They also appear in the statement of the \LR rule.
Note that some authors~\cite{Macdonald.1995,Sagan01}
prefer to write Yamanouchi words backwards and call them lattice permutations.

For a word $w$ we write $\abs{w}_x$ the number of occurrences of $x$ in $w$,
as computed by the \MC function \verb{count_mem x w}.
\begin{DEFN}
  A word $w := w_0,\dots,w_{l-1}$ over the integers is \emph{Yamanouchi} if for
  all $k, i \in \N$,
  \begin{equation}\label{eq.defYama}
    \abs{w_i,\dots,w_{l-1}}_k \geq \abs{w_i,\dots,w_{l-1}}_{k+1}\,.
  \end{equation}
\end{DEFN}
This is formalized by the following statement, where \verb|drop i s| 
removes the $i$-th first element of a sequence $s$:
\begin{lstlisting}
  forall i n, count_mem n (drop i s) >= count_mem n.+1 (drop i s)
\end{lstlisting}
The Yamanouchi words of length smaller than $4$ are exactly:
\begin{gather*}
  (), 0, 00, 10, 000, 100, 010, 210, \\
  0000, 1010, 1100, 0010, 0100, 1000, 0210, 2010, 2100, 3210
\end{gather*}
Applying \cref{eq.defYama} in the particular case $i=0$, we see that the
\emph{evaluation of $w$}, that is the sequence $(\abs{w}_i)$ where $i$ ranges
from $0$ to $\max(w)$ is decreasing a therefore partition of the size $l$ of
$w$.

For standard tableaux as well as for Yamanouchi words, we introduce types
(dependent pairs) for prescribed size, shape or evaluation, and endow then
with a canonical structure of \verb{fintype} structure by implementing
enumeration functions.

\section{Symmetric functions}
\label{sec:sf}
\subsection{Symmetric functions and Schur function}

We fix a positive integer $n$ and consider polynomials in the set of variables
$\alphX_n:=\{x_0,\dots,x_{n-1}\}$. The symmetric group $\SG_n$ on
$\{0,\dots,n-1\}$ acts on variables by permutation.
\begin{DEFN}
 A polynomial $f(x_0,\dots,x_{n-1})$ is symmetric if it is invariant by any
 permutation of the variables, that is for all permutation $\sigma\in\SG_n$,
 \begin{equation*}
   f^\sigma(\alphX) := f(x_{\sigma(0)}, x_{\sigma(1)}, \dots, x_{\sigma(n-1)}) = 
   f(x_{0}, x_{1}, \dots, x_{n-1})\,.
 \end{equation*}
\end{DEFN}
The sum of two symmetric polynomials is symmetric, as well as their product and 
the product of a symmetric polynomial by a scalar. Therefore, they form a
sub-algebra of the algebra of polynomials denoted $\sym(\alphX_n)$.

A natural basis of this algebra is given by the so-called monomial symmetric
polynomials. They are defined as the sum of the orbit of a monomial under the
action of the symmetric group. Since in any orbit there is only one monomial
whose exponents are sorted decreasingly, it makes sense to index the element
by this sequence of exponents. Removing the zeroes at the end of this sequence
gives a partition of length at most $n$. Therefore the bases of symmetric
polynomials are indexed by partitions of length at most $n$. Here is an
example of a monomial symmetric polynomial:
\begin{multline*}
  m_{(2,2,1)}(x_0,x_1,x_2,x_3) =
  x_0^2x_1^2x_2 + x_0^2x_1x_2^2 + x_0x_1^2x_2^2 \\
  + x_0^2x_1^2x_3 +  x_0^2x_2^2x_3 + x_1^2x_2^2x_3
  x_0^2x_1x_3^2 + x_0x_1^2x_3^2 + \\
  + x_0^2x_2x_3^2 + x_1^2x_2x_3^2 +
  x_0x_2^2x_3^2 + x_1x_2^2x_3^2\,.
\end{multline*}
It appears that this simple basis does not have that many interesting
properties. By contrast, one of the most important basis is the so-called
\emph{Schur polynomials}. They were first defined by Cauchy and Jacobi but
they are named in the honor of Schur who discovered their importance in the
representation theory. The original definition of Cauchy-Jacobi is the
following:
\begin{DEFN}[Cauchy-Jacobi's definition of Schur functions]
  Let $\lambda:=(\lambda_0\dots\lambda_{\ell-1})$ be a partition of length
  $\ell$ at most $n$. We complete $\lambda$ by setting $\lambda_i:=0$ whenever
  $i\geq\ell$. Then
  \begin{equation}
    s_{\lambda} (\alphX) :=
    \frac{1}{\Delta}\sum_{\sigma\in\SG_n}
    \sgn(\sigma)\ 
    \left(x_{\sigma(0)}^{\lambda_0 + n-1}
          x_{\sigma(1)}^{\lambda_1 + n-2}\cdots x_{\sigma(n-1)}^{\lambda_{n-1}}\right)
  \end{equation}
  where
  \begin{itemize}
  \item $\Delta:=\prod_{0\leq i<j<n} (x_i - x_j)$ is the Vandermonde
    determinant;
  \item $\sgn(\sigma)$ is $+1$ or $-1$ according to whether $\sigma$ is an
    even or odd permutation.
  \end{itemize}
\end{DEFN}
By (anti)-symmetry, the above quotient is a proper polynomial,
which is moreover symmetric. Furthermore it is easy to see that
\begin{equation}
s_{\lambda} (\alphX) =
   x_0^{\lambda_0}x_1^{\lambda_1}\cdots x_{n-1}^{\lambda_{n-1}} +  \cdots
\end{equation}
where the rest of the terms contains only monomials that are larger for the
reverse lexicographic order on the exponents. As a consequence, the Schur
polynomials are linearly independent and form a \emph{basis of the algebra of
  symmetric polynomials}.

This definition has an equivalent variant with a much more combinatorial flavor:
\begin{DEFN}[Combinatorial definition of Schur function]

  Let $\lambda$ be a partition of length
  at most $n$. For a tableau $t$ over the alphabet $\{0,\dots,n-1\}$, we
  denote $\alphX^t$ the product $\prod_{i\in t}x_i$ (equivalently it is the
  image of the row reading of $t$ under the morphism $i\mapsto x_i$ that
  sends noncommutative words to commutative monomials). Then
  \begin{equation}
    s_\lambda(\alphX) := \sum_{t\ \mid\ \shape(t) = \lambda}  \alphX^t\,.
  \end{equation}
  the sum being taken over all the tableaux of shape $\lambda$.
\end{DEFN}
For example, here is the computation of $s_{(2,1)}(x_0,x_1,x_2)$
\begin{gather*}
  \Yboxdim{9pt}
  \begin{array}[b]{c@{\ +\ }c@{\ +\ }c@{\ +\ }c@{\ +\ }c@{\ +\ }c@{\ +\ }l}
    \young(1,00) & \young(1,01) & \young(2,00) &
    \young(1,02)\ \young(2,01) & \young(2,11) & \young(2,02) & \young(2,12) \\[5mm]
    x_0^2x_1 & x_0x_1^2 & x_0^2x_2 & 2\,x_0x_1x_2 & x_1^2x_2 & x_0x_2^2 & x_1x_2^2\,.
  \end{array}
\end{gather*}
Though it is non trivial in general, one can check on this example that the
result is actually a symmetric polynomial. Indeed, one can subsume the previous
example by
\begin{equation*}
  s_{(2,1)} = 2\,m_{(1,1,1)} + m_{(2,1)}\,.
\end{equation*}
Here is a larger example
\begin{multline}\label{eq:example_schur}
  s_{3211} = 35\,m_{1111111} + 15\,m_{211111} + 6\,m_{22111} + \\
  2\,m_{2221} + 3\,m_{31111} + m_{3211}\,.
\end{multline}
In this expansion, the coefficient $6$ for $m_{22111}$, that is the
coefficient of $x_0^2x_1^2x_2x_3x_4$, is the number of tableaux of shape
$(3,2,1,1)$ whose row reading is a permutation of $0011234$. These tableaux
are:
\begin{equation*}
  \Yboxdim{10pt}\scriptscriptstyle
  \young(4,3,12,001)\quad\young(4,2,13,001)\quad\young(3,2,14,001)\quad
  \young(4,3,11,002)\quad\young(4,2,11,003)\quad\young(3,2,11,004)\,.
\end{equation*}
The fact that $s_{3211}$ is symmetric expresses for example that the coefficient
of $x_0^2x_1^2x_2x_3x_4$ is also the coefficient of $x_0x_1^2x_2x_3x_4^2$ as
one can check on the following list of tableaux
\begin{equation*}
  \Yboxdim{10pt}\scriptscriptstyle
  \young(4,3,24,011)\quad\young(4,3,14,012)\quad\young(4,3,12,014)\quad
  \young(4,2,13,014)\quad\young(4,2,14,013)\quad\young(3,2,14,014)
\end{equation*}
Again the reason why these two sets of tableaux are equinumerous is non-trivial.

Another very important point is that all these computations are independent from
the number of variables. More precisely, the equation of the preceding example
(\cref{eq:example_schur}) holds for any number of
variables, the difference being that $s_\lambda(\alphX_n)$ and
$m_\lambda(\alphX_n)$ vanish if $\ell(\lambda) > n$.

Now to formalize Schur polynomials, we need to introduce the commutative image of a word,
and the sum of the commutative image of a set a words of a given length
\verb{d}:
\begin{lstlisting}
Definition commword (w : seq 'I_n) : {mpoly R[n]} :=
  \prod_(i <- w) 'X_i.
Definition polylang d (s : {set d.-tuple 'I_n}) :=
  \sum_(w in s) commword w.
\end{lstlisting}
An important remark is that once its shape is fixed, one can recover any
(skew)-tableau from its row reading. We therefore define Schur function
$s_\lambda$ as the commutative image of the set of the row-reading of tableaux
of shape $\lambda$. Again, we define a decision procedure called
\verb|is_tableau_of_shape_reading| for characterizing such row-readings. Here is a
\verb|Prop| statement equivalent to this predicate:
\begin{lstlisting}
  exists tab, is_tableau tab /\ shape tab = sh /\ to_word tab = w.
\end{lstlisting}
Then, the definition of a Schur function is the following, where \verb|sh| is
a partition of \verb|d| ($\Sigma$-type \verb|intpartn d|):
\begin{lstlisting}
Definition tabwordshape (sh : intpartn d) :=
 [set t : d.-tuple 'I_n | is_tableau_of_shape_reading sh t].
Definition Schur d sh := polylang R (tabwordshape sh).
\end{lstlisting}

\subsection{The statement of the rule}

At last we come to the statement of the rule. It expresses the coefficients of
the product (the so-called structure constants) of the symmetric polynomials
in the basis of Schur polynomials. Let $A$ be a algebra over a field $\K$ with
a basis $B := (b_i)_{i\in I}$ indexed by a set $I$. The \emph{structure
  constants of $A$ in the basis $B$} are the coefficients $C_{i,j}^k$ with
$i,j,k\in I$ of the expansion of the product
$b_i b_j= \sum_{k\in I} C_{i,j}^k b_k$.

The \emph{\LR coefficients} $C_{\lambda, \mu}^{\nu}$ are the structure
constants of the algebra of symmetric polynomials in the basis of Schur
polynomials.

\begin{figure*}
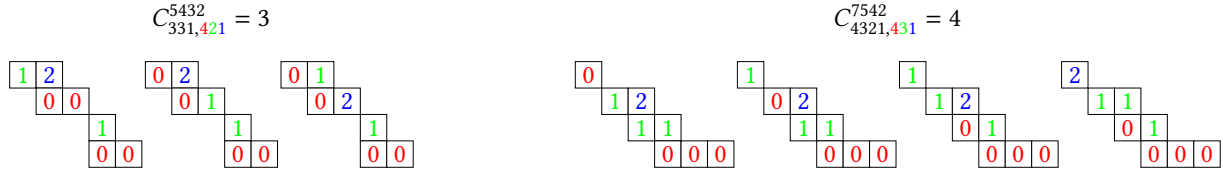

  \def\AA{\red 0} \def\AB{\grn 1} \def\AC{\blu 2} \def\AD{{\color{gray} 3}}
  \[
  \begin{array}{cccc}
    &C_{331,\red4\grn2\blu1}^{5432} = 3&
    &C_{4321,\red4\grn3\blu1}^{7542} = 4 \\[3mm]
    &\Yboxdim{10pt}\scriptscriptstyle
    \gyoung(\AB\AC,:;\AA\AA,:::;\AB,:::;\AA\AA)
    \gyoung(\AA\AC,:;\AA\AB,:::;\AB,:::;\AA\AA)
    \gyoung(\AA\AB,:;\AA\AC,:::;\AB,:::;\AA\AA)
    &\qquad\qquad
    &
    \Yboxdim{10pt}\scriptscriptstyle
    \gyoung(:;\AA,::;\AB\AC,:::;\AB\AB,::::;\AA\AA\AA)
    \gyoung(:;\AB,::;\AA\AC,:::;\AB\AB,::::;\AA\AA\AA)
    \gyoung(:;\AB,::;\AB\AC,:::;\AA\AB,::::;\AA\AA\AA)
    \gyoung(:;\AC,::;\AB\AB,:::;\AA\AB,::::;\AA\AA\AA)
  \end{array}
  \]
  \caption{Two examples of Littlewood-Richardson coefficients}
  \label{fig:LR-example}
\end{figure*}

The \LR rule states that the coefficients are non-negative integers, and
gives a way to compute them:
\begin{theorem}[Littlewood-Richardson rule]\label{theo:LR-rule}
  $C_{\lambda, \mu}^{\nu}$ is the number of (skew) tableaux of shape the
  difference $\nu/\lambda$, whose row reading is a Yamanouchi word of
  evaluation $\mu$.
\end{theorem}
Recall that the evaluation condition means that, for all $i$, the word $w$
contains exactly $\mu_i$ many occurrences of the letter
$i$. See~\cref{fig:LR-example} for some examples.
Since once the shape is fixed, one can recover a tableau from its row
reading, the statement can be rephrased:
\begin{theorem}[Littlewood-Richardson rule]
  $C_{\lambda, \mu}^{\nu}$ is the number of Yamanouchi words of evaluation
  $\mu$ that are the row reading of a (skew) tableaux of shape the difference
  $\nu/\lambda$.
\end{theorem}

Recall that in \MC, speaking of the cardinality of a set, requires its
elements to inhabit a finite type. We thus provide a type named \verb{yameval}
for Yamanouchi words of prescribed evaluation $\mu$, with a canonical
structure of {\verb|fintype|.} Recall also that partitions of an integer
\verb+d+ are encoded by terms of type \verb|intpartn d|. The following
definitions and reflection lemma formalize the property of being the
row-reading of a skew tableau of shape $\var{P}/\var{P1}$:
\begin{lstlisting}
Definition is_skew_reshape_tableau P P1 (w : seq nat) :=
  is_skew_tableau P1 (skew_reshape P1 P w).
Lemma is_skew_reshape_tableauP P P1 (w : seq nat) :
  size w = sumn (diff_shape P1 P) -> reflect
    (exists tab, is_skew_tableau P1 tab
                 /\ shape tab = diff_shape P1 P
                 /\ to_word tab = w)
    (is_skew_reshape_tableau P P1 w).
\end{lstlisting}
Now the \LR coefficient $C_{\var{P1},\var{P2}}^{\var{P}}$ is the cardinality
of the set of Yamanouchi words of evaluation~\var{P2} that are row reading of
tableaux of shape $\var{P}/\var{P1}$:
\begin{lstlisting}
Definition LRyam_set P P1 P2 :=
  [set x : yameval P2 | is_skew_reshape_tableau P P1 x].
Definition LRyam_coeff P P1 P2 := #|LRyam_set P P1 P2|.
\end{lstlisting}
With all these definitions, the \LR rule statement is:
\begin{lstlisting}
  Schur P1 * Schur P2 =
  \sum_(P : intpartn (d1 + d2) | included P1 P)
    Schur P *+ LRyam_coeff P P1 P2.
\end{lstlisting}
We will discuss some implementations of the rule in~\cref{sec:implem}.

\section{Character theory of the symmetric groups}
\label{sec:SymGroup}

In this section we give a broad overview of formalized results concerning the
character theory of the symmetric groups. We first recall quickly some very
basic definition of character theory of a group.

For a finite group $G$, two elements $a, b\in G$ are \emph{conjugate} if there
is a $g\in G$ such that $a = gag^{-1}$. This is an equivalence relation whose
equivalence classes are called \emph{conjugacy classes}. Functions from $G$ to
the field of complex number $\mathbb{C}$ that are constant for members of the
same conjugacy class are called \emph{class functions}. A group representation
of degree $d$ is a map $\rho$ that associates to any element $g\in G$ an
invertible $d\times d$ matrix $\rho(g)$ such that $\rho(gh) = \rho(g)\rho(h)$,
otherwise said, a group homomorphism $g\to GL_d(\mathbb{C})$. The character of
the group representation $\rho$ is the class function that associates to each
group element the trace of the corresponding matrix. All these notion were
already available in \MC.
\medskip

All our results concerns specifically the the symmetric group $\SG_n$ of
permutations of the set $\{0, \dots, n-1\}$ (to follow \MC conventions).  Two
permutations $\sigma$ and $\mu$ are conjugate if and only if, the lengths of
their cycles in the canonical cycle decomposition agree. The collection of the
lengths of the cycle of $\sigma$ (including fixed points) is a partition of
$n$ called the cycle type of $\sigma$. This is the main first result (where
cycle type is a map from \verb|'S_n| to \verb|intpartn n|):
\begin{lstlisting}
(exists g, t = g * s * g ^-1) <-> (cycle_type s = cycle_type t).
\end{lstlisting}
By comparing the size of the conjugacy classes with the coefficients of the
expansion of the powersum basis in the complete basis of symmetric functions,
Frobenius associated to any class function of the symmetric group a symmetric
function called its Frobenius characteristic. This map happens to be an
isomorphism for most interesting structures (product, coproduct, scalar
product, Kronecker product).

In our formalization, the map is denoted \verb|FChar|. The following 
statement asserts that it is a ring morphism (sending the induction
product of class function to the usual product of polynomials).
Here \verb|'CF('SG_m)| denote the type of class function of the
symmetric groups $\SG_m$ and \verb|'Ind['SG_(m + n)]| is the induction map of
class function from $\SG_{m}\times\SG_{n}$ to $\SG_{m + n}$:
\begin{lstlisting}
Theorem Fchar_ind_morph (f : 'CF('SG_m)) (g : 'CF('SG_n)) :
  Fchar ('Ind['SG_(m + n)] (f \o^ g)) = Fchar f * Fchar g.
\end{lstlisting}
There are moreover two scalar products respectively on class functions and
symmetric functions for which the Frobenius characteristic is an isometry:
\begin{lstlisting}
Theorem Fchar_isometry f g : '[Fchar f | Fchar g] = '[f, g].
\end{lstlisting}
This allows to translate many statement about symmetric functions to
characters of the symmetric groups. Here is for example the translation of the
\LR rule in this language (\verb|la| and \verb|mu| are partition of
respectively $c$ and $d$ and \verb|'irrSG[la]| denotes the irreducible
character associated to \verb|la|):
\begin{lstlisting}
'Ind['SG_(c + d)] ('irrSG[la] \o^ 'irrSG[mu]) =
  \sum_(nu : 'P_(c + d) | included la nu)
          'irrSG[nu] *+ LRyam_coeff la mu nu.
\end{lstlisting}
Using a backtracking search similar to the one used in the \LR rule, we also
provide an implementation of the Murnaghan-Nakayama rule to compute the value
of the irreducible characters. This is ensured by the following theorem:
\begin{lstlisting}
Theorem Murnaghan_Nakayama_char n la (sigma : 'S_n) :
  'irrSG[la] sigma = MN_coeff la (cycle_type sigma).
\end{lstlisting}
This allow to actually compute the value in Coq/Rocq. We can check the value
given as example on Wikipedia:
\begin{lstlisting}
Example Wikipedia_Murnaghan_Nakayama :
  let p521 := @IntPartN 8 [:: 5; 2; 1]%N is_true_true in
  let p3311 := @IntPartN 8 [:: 3; 3; 1; 1]%N is_true_true
     in 'irrSG[p521] (permCT p3311) = - 2%:~R.
Proof. by rewrite /= Murnaghan_NakayamaCT. Qed.
\end{lstlisting}
See \cref{sec:implem} for some remarks about the implementations.

\section{Overview of the library}
\label{sec:overview}
The library was adopted by the \texttt{math-comp} organization in October 2021
and is accessible on \url{https://github.com/math-comp/Coq-Combi}
together with its documentation
\url{https://math-comp.github.io/combi/1.0.0/toc.html} and the precise
requirement to be able to compile it.

The contributed library can be divided into a first combinatorial part,
dealing with various combinatorial objects, and two algebraic parts, dealing
respectively with symmetric functions and the representation theory of the
symmetric group. We moreover provide a handful of side results that are of
general interest, related to the present theories. See \cref{appendix:toc} and
\cref{fig.depend} for the detailed table of contents of the library.

\paragraph{Combinatorial objects}

\begin{itemize}
\item integer partitions, together with Young's and dominance
  lattices, skew partitions, horizontal, vertical and ribbon border strips;
\item integer composition, bijection with subsets and refinement lattice;
\item subsequences, integer vectors;
\item standard tableaux, skew tableaux, Yamanouchi word;
\item permutations, standard words, and the standardization map;
\item binary trees, Dyck words, Tamari lattices and Catalan numbers;
\item set partitions and refinement lattices;
\end{itemize}

\paragraph{Theory of symmetric functions}
\begin{itemize}
\item Schur function and Kostka numbers and the equivalence of the
  combinatorial and algebraic (Jacobi) definition of Schur polynomials;
\item the classical bases, Newton formulas and various basis changes;
\item the scalar product and the Cauchy formula;
\item The Littlewood-Richardson rule using Schützenberger approach, itself including:
\item the \RS correspondence;
\item the construction of the plactic monoid using Greene invariants;
\item the Littlewood-Richardson and Pieri rules;
\item the Murnaghan-Nakayama rule for converting power sums to Schur
  functions.
\end{itemize}

\paragraph{The symmetric Group and its character theory}

This development was done with contributions from Thibaut Benjamin. Most of the
results about symmetric functions can be rephrased in terms of character
theory of the symmetric group using Frobenius isomorphism. Here is an overview
of the contents of this part of the library:
\begin{itemize}
 \item the Coxeter presentation of the symmetric group, Matsumoto theorem stating
   that the symmetric group $\SG_n$ is generated by the elementary
   transpositions $s_i$ exchanging $i$ and $i+1$, subject to the braid
   relations;
\item the Coxeter length and the inversion set;
\item the weak permutohedron lattice;
\item cycle types for permutations: a permutation admits a unique
  decomposition (up to the order) as a product of cycle with disjoint
  supports;
\item the tower structure and the restriction and induction formulas for class
  indicator;
\item the structure of the centralizer of a permutation and the cardinality of
  the conjugacy classes of the symmetric group;
\item Young's characters and Young's Rule;
\item the theory of Frobenius characteristic and Frobenius character formula;
\item the Murnaghan-Nakayama and  Littlewood-Richardson rules in terms of
  group characters.
\end{itemize}

\paragraph{Various side results of wider interest}
\label{par.side-results}
\begin{itemize}
\item the Erdös Szekeres theorem about increasing and decreasing
  subsequences from Greene's invariants theorem~\cite{Erdos87}. This is a
  simple consequence of the \RS correspondence. The theorem appears in
  ``Formalizing 100 Theorems''\cite{Theorem100};
\item the factorization of the Vandermonde determinant as the product of
  differences (used for example in Lagrange interpolation);
\item the Hook-Length Formula giving the number of standard Young tableaux, closely
  following the probabilistic proof of~\cite{Greene79} (done in collaboration
  with Christine Paulin and Olivier Stietel).
\end{itemize}

\section{Outline of the proof : the Robinson-Schensted algorithm and the
  plactic monoid}
\label{sec:LRproof}
The formalized proof roughly follows the Schützenberger argument as presented
presented in~\cite{Lothaire.2002-Plact}. It is based on an in-depth study of a
classical algorithm due to Schensted, which computes the length of a longest
non-decreasing subsequence of a word. In particular, the central argument is a
description of the output of the algorithm on the concatenation of two words,
knowing the output on those. Therefore -- a typical feature of algebraic
combinatorics -- this is a proof of an algebraic identity, based on the
understanding of the behavior of an algorithm. The main steps of the proofs are:
  \begin{enumerate}
  \item increasing subsequences and Schensted's algorithms;
  \item Robinson-Schensted correspondence: a bijection;
  \item Greene invariants: computing the maximum sum of the length of $k$ disjoint
    non-decreassing subsequences;
  \item Knuth relations, the plactic monoid;
  \item Greene invariants are plactic invariants: equivalence between RS and plactic;
  \item standardization; symmetry of RS;
  \item lifting to non commutative polynomials: free quasi-symmetric function
    and shuffle product;
  \item noncommutative lifting of the LR-rule: the free/tableau LR-rule;
  \item back to Yamanouchi words: a final bijection.
  \end{enumerate}
  In this section, we briefly present a few selected steps. We start with the main
  algorithmic ingredient, namely the \RS correspondence.
  
\subsection{Schensted algorithm and the \RS bijection}

\begin{figure*}[ht]
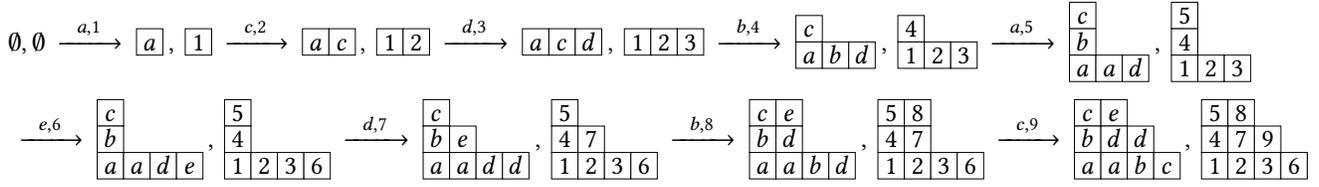

  \Yboxdim{10pt}
  \def\ar#1{\ \xrightarrow{\ \ #1\ \ }\ }
\begin{multline*}
\emptyset, \emptyset
\ar{a,1}
\young(a)\,,\ \young(1)
\ar{c,2}
\young(ac)\,,\ \young(12)
\ar{d,3}
\young(acd)\,,\ \young(123)
\ar{b,4}
\young(c,abd)\,,\ \young(4,123)
\ar{a,5}
\young(c,b,aad)\,,\ \young(5,4,123)
\\
\quad\ar{e,6}
\young(c,b,aade)\,,\ \young(5,4,1236)
\ar{d,7}
\young(c,be,aadd)\,,\ \young(5,47,1236)
\ar{b,8}
\young(ce,bd,aabd)\,,\ \young(58,47,1236)
\ar{c,9}
\young(ce,bdd,aabc)\,,\ \young(58,479,1236)
\end{multline*}
\caption{\RS algorithm applied to the word $acdbaedbc$}
\label{fig:exampleRS}
\end{figure*}

\begin{DEFN}
  Let $w = w_0\dots w_{l-1}$ a finite sequence. A \emph{subsequence} of $w$ is
  a sequence $v=v_0\dots v_{r-1}$ such that there exists a set of integers
  $S:=\{i_0<i_1<\dots<i_{r-1}\}$ verifying $0\leq i_0$ and $i_{r-1}< l$ and
  for all $k<r$ then $v_k= w_{i_k}$. The set $S$ is called the support of the
  subsequence.
\end{DEFN}
Subsequence are  already defined in \MC by the following boolean predicate:
\begin{lstlisting}
Fixpoint subseq s1 s2 :=
  if s2 is y :: s2' then
    if s1 is x :: s1' then subseq (if x == y then s1' else s1) s2' else true
  else s1 == [::].
\end{lstlisting}
Now let us consider the following problem:
\begin{problem}
  Given a finite sequence $w$ over a totally ordered set, compute the maximum
  length of a non-decreasing subsequence.
\end{problem}
We stress that the problem is to compute the length, but not necessarily an
actual longest non-decreasing subsequence.
The problem is solved by a kind of dynamic programming algorithm due to
Schensted~\cite{Schensted61}:
\begin{ALGO}\label{algo:Schensted}
  Start with an empty row $r$, insert the letters $l$ of the word one by one
  from left to right by the following rule:
  \begin{itemize}
  \item replace the first letter strictly larger than $l$ by $l$;
  \item append $l$ to $r$ if there is no such letter.
  \end{itemize}
\end{ALGO}
This is easily implemented in \Coq:
\begin{lstlisting}
Fixpoint insrow r l : seq T :=
  if r is l0 :: r then
    if (l < l0)%Ord then l :: r
    else l0 :: (insrow r l)
  else [:: l].
Fixpoint Sch_rev w :=
  if w is l0 :: w' then ins (Sch_rev w') l0 else [::].
Definition Sch w := Sch_rev (rev w).
\end{lstlisting}
Then Schensted showed the following invariant:
\begin{theorem}
  The $i$-th entry $r[i-1]$ of the row $r$ is the smallest letter that ends a
  non-decreasing subsequence of length~$i$. As a consequence, the length of $r$ is
  the maximum length of a non-decreasing subsequence.
\end{theorem}
Here is the formalized statement, where \verb|subseqs w| is a type for the subsequences of sequences \verb|w|:
\begin{lstlisting}
Theorem Sch_max_size (w : seq T) :
  size (Sch w) = \max_(s : subseqs w | sorted s) size s.
\end{lstlisting}
Note that in the first case of the Schensted algorithm the first letter $k$
in row $r$ strictly larger than $l$ is replaced by $l$. We say that $k$ is
\emph{bumped} by $l$. We can then insert the bumped letters using the same algorithm
in a new row $r_1$ (defining $r_0:=r$). If inserting letters in $r_1$ bumps
some letters, they are inserted in $r_2$, an so on. By stacking the rows, one
gets a proper tableau called \emph{the insertion tableau of $w$} usually
denoted by $P(w)$. If one records in another tableau $Q(w)$ (called the
\emph{recording tableau}) the order in which the boxes were added, one gets a
pair of tableaux of the same shape, the second one being standard. This is the
\RS correspondence, exemplified in~\cref{fig:exampleRS}.
The main point of the recording tableau is that it allows to reverse the whole
construction and recover $w$ from the pair $(P(w), Q(w))$:
\begin{theorem}
  The map $w \mapsto (P(w), Q(w))$ is a bijection between word on $\alphA$ and
  pairs $(P, Q)$ where
  \begin{itemize}
  \item $P$ is a tableau on $\alphA$;
  \item $Q$ is a standard tableau;
  \item $P$ and $Q$ have the same shape.
  \end{itemize}
\end{theorem}
Here is the formalized statement:
\begin{lstlisting}
Definition is_RStabpair pair := let: (P, Q) := pair in
  is_tableau P /\ is_stdtab Q /\ (shape P == shape Q).
Structure rstabpair :=  (* Sigma type *)
  RSTabPair { pqpair :> TabPair; _ : is_RStabpair pqpair }.
Theorem bijRStab : bijective (RStab : seq T -> rstabpair).
\end{lstlisting}
A natural question at this point is: what does the \RS algorithm compute? The
question was answered by Greene~\cite{Greene74}.  Let $w$ be a word. We
consider all the tuples of non decreasing subsequences $(v_1, \dots, v_l)$ whose
associated supports $(S_1,\dots,S_k)$ are disjoint. For examples, the green red
and blue subsequences below are respectively associated to the disjoint sets
$\{0,1,3,4\}$, $\{2,5,6,7,9\}$, $\{8,10,12\}$.
  \[
  \begin{array}{c@{}c@{}c@{}c@{}c@{}c@{}c@{}c@{}c@{}c@{}c@{}c@{}c}
    \grn{a}&\grn{b}&\red{a}&\grn{b}&\grn{c}&\red a&\red b&\red b&\blu{a}&\red d&\blu{b}&a&\blu{b}
  \end{array}
\]
Then Greene showed that:
\begin{theorem}
  For any word $w$, and integer $k$, the sum of the length of the first $k$
  rows of $P(w)$ is the maximum sum of the length of $k$ disjoint support
  non-decreasing subsequences of $w$.
\end{theorem}
Note that an analogue statement exists for strictly decreasing subsequences. Its
formalization relies on the following definition, where \verb|cover S| is the
union of the set of sets \verb|S|, \verb|trivIset S| means that the
elements of \verb|S| are pairwise disjoint and \verb|extract w s| is
the subsequence extracted from \verb|w| with support \verb|s|:
\begin{lstlisting}
Definition ksupp k (S : {set {set 'I_N}}) : bool :=
  #|S| <= k && trivIset S
  && [forall (s | s \in S), sorted (extract w s)].
\end{lstlisting}
Here as well we follow the small scale reflection approach, and provide an executable definition: \verb+&&+ is an infix notation for boolean conjunctions, and the bracketed quantification, a bounded one. Indeed, the quantified variable \verb|s| ranges over the finite type \verb|{set 'I_N}| of bounded integers, filtered by membership to \verb+S+.  With this definition, Green's theorem can be written as:

\begin{lstlisting}
Definition greenRow := \max_(S | ksupp k S) #|cover S|.
Theorem greenRow_RS k w :
  greenRow w k = \sum_(l <- take k (shape (RS w))) l.
\end{lstlisting}
where \verb|RS(w)| is the Coq notation for $P(w)$. This ingredient is the crux of Knuth's theorem~\cite{Knuth70}. Recall that
a congruence $\equiv$ on the words on $\alphA$ is an equivalence relation
which is preserved by left and right concatenation (\ie $u\equiv v$ implies
$uw \equiv vw$ and $wu \equiv wv$). In this case, the quotient of the free
monoid $\alphA^*$ by $\equiv$ is itself endowed with a structure of monoid.
\begin{theorem}\label{theo:plact}
  Consider the smallest congruence $\placteq$ that contains $bca \placteq
  bac$ whenever $a < b \leq c$ and $acb \placteq cab$ whenever $a \leq b < c$
  where $a,b,c$ are any letter in $\alphA$. Then for any two words $v, w$ the
  equivalence $v \placteq w$ holds if and only if $P(v)=P(w)$.
\end{theorem}
Theorem~\ref{theo:plact} is formalized as a boolean equality:
\begin{lstlisting}
Theorem plactic_RS u v : (u =Pl v) = (RS u == RS v).
\end{lstlisting}
The full formal definition of \verb|=Pl| is however too long to be displayed
here, involving a search algorithm to compute a connected component in a
finite graph. It's specification are written thanks to statements such as
\begin{lstlisting}
  exists a b c, a <= b < c /\ u = [:: a; c; b] /\ v = [::c; a; b]
\end{lstlisting}
The reverse direction of Knuth's theorem is relatively easy to prove. By
following step by step what happens in the execution of the \RS algorithm, one
can show that $w$ is plactic equivalent to the row reading of $P(w)$.
\begin{lstlisting}
  Theorem congr_RS w : w =Pl (to_word (RS w)).
  Corollary Sch_plact u v : RS u == RS v -> u =Pl v.
\end{lstlisting}
The direct implication uses Green's theorem to show that two plactic equivalent
words have tableaux of the same shapes:
\begin{lstlisting}
  u =Pl v -> shape (RS u) = shape (RS v).
\end{lstlisting}
The actual equality of \verb|RS u| and \verb|RS v| follows by some
specifically tailored induction following the removing of the last occurrence
of the largest letter of a word.
\medskip

The hard part of the proof of Green's theorem amounts to proving that if two
words differ by a plactic rewriting, their $k$-Green numbers (maximum sum of
the length of $k$ disjoint support non-decreasing subsequences of $w$) are
equal. We prove the latter by using the following predicate that asserts that
for any $k$-support of \verb|u1| there exists a $k$-support of \verb|u2| with
the same cover size.
\begin{lstlisting}
Definition ksupp_inj k (u1 : seq T1) (u2 : seq T2) :=
  forall s1, s1 \is a k.-supp[in_tuple u1] ->
    exists s2, (#|cover s1| == #|cover s2|)
            /\ (s2 \is a k.-supp[in_tuple u2]).
Lemma ksuppRow_inj_plact1i : v2 \in plact1i v1 ->
  ksupp_inj k (u ++ v1 ++ w) (u ++ v2 ++ w)
Corollary GreeneRow_leq_plact1i : v2 \in plact1i v1 ->
  GreeneRow (u ++ v1 ++ w) k <= GreeneRow (u ++ v2 ++ w) k.
\end{lstlisting}
This is proved by some surgery in the $k$-support. This part of the proof was
actually quite hard to mechanize. First of all, to avoid duplicating arguments
in the proof, we abstracted the various cases using \Coq's modules. Here is
some details on one case (module \verb|SetContainingBothLeft|): we denote
$x:=ubacv$ and consider a given $k$-support $P$ containing both $a$ and
$c$. We want to deal at once with the case where $a < b \leq c$ and
$c < b \leq a$. So we encaspulate the hypothesis in a record \verb|hypRabc|
(where $R:=\leq$ in the first case and $R:=>$) in the second case):
\begin{lstlisting}
Record hypRabc R a b c := HypRabc {
   Rtrans : transitive R;
   Hbc : R b c;  Hba : ~~ R b a;
   Hxba : forall l, R l a -> R l b;  Hbax : forall l, R b l -> R a l
}.
\end{lstlisting}
Then in the case where \verb|b| is not in \verb|cover P|, replacing \verb|a|
by \verb|b| gives another $k$-support. On the other hand if \verb|b| belongs
to \verb|cover P| one has to exchange some right part of $k$-support using
surgery such as
\begin{lstlisting}
      (S :&: [set j : 'I_(size x) | j <= posa])
  :|: (T :&: [set j : 'I_(size x) | j > posc]).
\end{lstlisting}
where \verb|posa| and \verb|posc| are the respective position of \verb|a| and
\verb|c| in $x:=ubacv$ and \verb|:&:| and \verb':|:' are \MC notations for the
 intersection and union of finite sets. The careful reader may have observed that
this definition depends on the underlying word \verb|x|. This admittedly might have been an unfortunate design choice, as it incurs a lot of hurdles. Anyway, we can state and prove this this way theorems such as:
\begin{lstlisting}
Let x := u ++ [:: b; a; c] ++ v.
Let y := u ++ [:: b; c; a] ++ v.
Hypothesis HRabc : hypRabc R a b c.
Variable P : {set {set 'I_(size x)}}.
Hypothesis Px : P \is a k.-supp[R, in_tuple x].
Theorem exists_Qy : exists Q : {set {set 'I_(size y)}},
  #|cover Q| = #|cover P| /\ Q \is a k.-supp[R, in_tuple y].
\end{lstlisting}
It claims that for any $k$-support of $x:=ubacv$, if the rewriting to
$y:=ubcav$ use a valid plactic relation, there exists a $k$-support for $y$
with the same cover size. This ensure that the $k$-green number of $y$ is at
least the one of $x$. Similar statement deals with all plactic rule,
completing the proof of Green's and Knuth's theorems.

\subsection{Noncommutative lifting of Schur function}
\label{subsec:noncomm}
In the previous section, we have already seen some interaction between
combinatorics and algebra. Indeed Knuth's theorem basically says that the \RS
algorithm computes a normal form for the plactic monoid~\cite{SchutzLR}. Note
that this is not an isolated case. For example, the well known algorithm of
binary search tree insertion also computes a normal form of a quotient of the
free monoid, called the Sylvester monoid~\cite{Hivert2005}. The aim of this
section is to present yet another example of such interaction, namely the
central argument of the proof of the \LR rule.
\smallskip

The main idea is to lift Schur functions at a noncommutative level, where
multiplication of polynomials becomes concatenation of languages. Recall that
all considered languages are homogeneous and can thus be encoded by terms of
type \verb|homlang d := {set d.-tuple 'I_n}|. Recall as well that the polynomial \verb|polylang s| is defined as the sum of the commutative image of the set \verb|s|. The
concatenation of two languages can be defined as follows, for \verb|d1, d2 : nat|:
\begin{lstlisting}
Definition catlang (s1 : homlang d1) (s2 : homlang d2) :
  homlang (d1 + d2) :=
  [set cat_tuple w1 w2 | w1 in s1, w2 in s2].
\end{lstlisting}
The \verb|polylang| function sends concatenations of languages to commutative products of
polynomials as stated by
\begin{lstlisting}
  polylang s1 * polylang s2 = polylang (catlang s1 s2).
\end{lstlisting}
One way to prove \LR would be to concatenate two languages whose commutative
image is the Schur function. By the very definition we gave of Schur functions,
such a language is \verb|tabwordshape sh|. However it is not the only
one. Indeed, thanks to the \RS correspondance, we have as many elements in a
plactic class as standard tableaux of its shape. More precisely, if \verb|Q| is a fixed
standard tableau of shape \verb|sh|, then for any tableau \verb|P| of the same
shape, the word $w$ that corresponds to the pair \verb|(P, Q)| inserts to
\verb|P| and thus has the same commutative image as the row reading of
\verb|P|. Therefore the Schur function is nothing but the commutative image of
the set of words $w$ whose recording tableau is \verb|Q|. We call this set a
free Schur function.
\begin{lstlisting}
Definition freeSchur (Q : stdtabn d) : homlang d  :=
  [set t : d.-tuple 'I_n | (RStabmap t).2 == Q].
\end{lstlisting}
The above reasoning shows that:
\begin{lstlisting}
  Schur (shape_deg Q) = polylang R (freeSchur Q).
\end{lstlisting}
Note that there are as many non-commutative lifting of the Schur function
associated to $\lambda$ as standard tableaux of shape $\lambda$.  Remark that
the free Schur languages form a set partition of the set of all words. The
crucial fact is that this set partition is stable by concatenation meaning that
the concatenation of two parts is a union of parts. So there exists a certain
set \verb|LRsupport Q1 Q2| such that
\begin{lstlisting}
  catlang (freeSchur Q1) (freeSchur Q2) =
     \bigcup_(Q in LRsupport Q1 Q2) freeSchur Q.
\end{lstlisting}
The definition of \verb|LRsupport| involves shuffling words and is too long
to be described here. If \verb|sh| is a partition of \verb|d1 + d2|, we can define
\begin{lstlisting}
Definition LRtab_set Q1 Q2 sh :=
  [set Q in LRsupport Q1 Q2 | shape Q == sh].
\end{lstlisting}
Then the cardinality of \verb|LRtab_set Q1 Q2 sh| is the \LR coefficient
associated to the shapes of \verb|Q1|, \verb|Q2| and the partition
\verb|sh|. In particular it doesn't depend on the particular choice of
\verb|Q1| and \verb|Q2| but only their shape. This is proved by an explicit
bijection. To complete the proof of the \LR rule, the last step consists in defining a bijection between the set \verb|LRsupport Q1 Q2| and the \LR
tableaux. Defining and proving such a bijection takes 700~extra lines of
code\dots

\section{Implementations}
\label{sec:implem}

Recall that the definition of the \LR coefficients we gave is the cardinality
of a set of Yamanouchi words verifying some condition. Our formal definition is however proof-oriented, rather than algorithmic: although effective, it uses constants whose reduction has been blocked by design, as is customary in \MC. Controlling redudction is indeed key to the robustness of notations, to the behavior of type inference, of term comparison, etc. As a consequence, \Coq is not able to compute a coefficient from this
definition. This is however very easily circumvented using filtering of enumeration
sequences as in:
\begin{lstlisting}
[seq x <- enum_yameval P2 | is_skew_reshape_tableau P P1 x]
\end{lstlisting}
Then the row readings of the four tableaux of the second example of
\cref{fig:LR-example} are recovered by:
\begin{lstlisting}
Eval compute in (LRyam_enum  [:: 4; 3; 2; 1] [:: 4; 3; 1] [:: 7; 5; 4; 2]).
= [:: [:: 0; 1; 2; 1; 1; 0; 0; 0];
      [:: 1; 0; 2; 1; 1; 0; 0; 0];
      [:: 1; 1; 2; 0; 1; 0; 0; 0];
      [:: 2; 1; 1; 0; 1; 0; 0; 0]] : seq (seq nat)
\end{lstlisting}
Note that this is not an efficient way of computing them. One cannot get
coefficients larger than 10 in a reasonable time using this definition. A more
efficient way is obtained using a backtracking search. We implemented one
directly in \Coq under the name \verb|LRcoeff| and proved that it computes the
same values as \verb|LRyam_coeff|. It is already relatively
efficient:
\begin{lstlisting}
Eval vm_compute in LRcoeff
  [:: 9; 6; 5; 4; 3; 2; 1] [:: 7; 6; 5; 5; 4; 3; 2; 1]
  [:: 11; 10; 9; 8; 7; 5; 4; 3; 3; 2; 1].
= 81672 : nat
\end{lstlisting}
This was computed by \Coq in $7.7s$. The extracted OCaml program computes the
same number in only $0.29s$. An optimized C implementation such as
\verb|lrcalc|~\cite{lrcalc} which is the one used by \Sage computes the same
number in $0.01s$. This overhead is expected, in particular due to the use of
mutable arrays instead of list. These timings are very small and seem to
indicate that there is no need to optimize things. However, in many practical
applications one needs not only to compute a single coefficient but rather to
expand a full product of two Schur functions. This amounts computing an
exponentially growing number of coefficients. So the computation is performed
by the same kind of backtracking that we are using, but with a unique
backtrack search to compute all the different coefficients (\ie with a varying
outer shape of the \LR tableau).

To conclude this section, we mention that several other computational rules
are implemented in our contributed library, notably the Murnaghan-Nakayama
rule for computing character values of the symmetric group. We are not aware
of any other formally verified implementations thereof.

\section{Concluding remarks}

In this paper we have presented the first (as far as we now) machine checked
proof as well as implementation of the \LR rule. Having these is very relevant since
there seems to exists no fundamentally better way to compute those numbers. We
would like to underline a few facts we learned during this development. They
may not come as a surprise for specialists, but we think they may be very
helpful, should the algebraic combinatorics community start to develop proofs
on a large scale as it has done for computations.

\smallskip In our point of view, the most important outcome is that it is
feasible!  Though we had a strong experience in the development of computer
algebra systems, we had none in participating in large formal proofs. We
neither had any formal training in \Coq or any proof assistant. The outcome is
currently more than 35k lines of code whose development spans several years.

Another very strong conclusion is that proving things in algebra is much
easier than in combinatorics. Indeed, algebraic structures are very polished
and lots of thinking has been done about cleaning up their
definition. Admittedly much work had already been done by the \MC project, but
the fact that the proof usually boils down to defining the right structure
followed by a series of rewriting helps a lot. On the other hand, for
combinatorics very little if anything can be shared from one development to an
other. This explains the very large number of functions in a computational
library (see below some figures). Moreover there are many places where a whole
theory is needed to make a choice where any would do. For example, given two
maps \verb|f:A -> nat| and \verb|g:B -> nat| whose cardinalities of preimages
of any integer agree, there is clearly a bijection \verb|r: A -> B| such that
\verb|g = f o r|.  Constructing such a bijection is non trivial (see the
\verb|Combi/fibered_set.v| file).

Except for definitions (see below), the proofs followed closely the
paper ones, with one big exception. There were no formal power series in
\MC. Recall that all the algebraic structures in \MC assume that equality is
decidable. This is a problem with power series as they are infinite objects.
So the authors decided to use a bijective approach when power series are
commonly used. This includes proofs of many counting formulas as the size of
conjugacy classes and the change of basis for symmetric functions such as
Newton formulas. After the main development, the authors wrote a library for
formal power series but the present development still needs to be backported.

\smallskip
\MC proved to be a very good platform for such a development. The choice was
quite obvious for the second part of the library which deals with the
character theory of the symmetric group as all the basic group and
representation theory was formalized during the proof of the Feit-Thomson odd
order theorem~\cite{GonthierFT}. However, this only came at a late stage of
the development, when we realized that all the statement we had expressed in
terms of polynomials could be relatively easily translated in term of
representation theory thanks to Frobenius isomorphism.

Most importantly the usage of small scale reflection is extremely effective in
combinatorics. Having the system decide by itself many particular case has
proven very effective to reduce the size of the proof allowing the developer
to focus on the important cases. Also, combinatorial definitions are very easy
to get wrong so that being able to test them is quite useful. Finally, since
decisions are programmed in \Coq itself, there is no need for high-powered
tactics. Except for a of handful use of \verb|repeat| no tactical programming
was used. The only advanced tactics we find useful are \verb|ring| and
\verb|lia| which weren't available in \MC at the time of writing the library,
so we had to hold the hand of the system during computations.

There is however a design choice of \MC that was unfortunate for us. \MC was
a little bit too geared for finite sets. The fact that sets can only be
defined on a finite type forces to use dependent types such as
\verb|d.-tuple T| instead of \verb|seq T| and this complicate defining
bijections and algorithms. Fortunately, this restriction is currently being
lifted with the development of \MC analysis.

\smallskip
One of the outcomes for the combinatorics community concerns definitions. It is
striking that the formal definition often departs largely from the paper
ones. Moreover, in a few cases, formalizing them makes us realize that it is
very common in the paper that a definition for a notion is given but clearly
the proofs use a different definition which is not given! Let us give two
examples that are simple enough to be explained in a few lines. The first
example is the notion of ribbon shape. Here is a definition taken from one of
the authors' paper, to avoid giving the impression of blaming a colleague:
\begin{quotation}
  A ribbon is a connected skew shape that doesn't contain any $2\times 2$
  square of boxes.
\end{quotation}
Here is an example and two counter examples (the second is not connected and
the third contains a $2\times 2$ square of boxes:
\[\Yboxdim{6pt}
  \gyoung(;;,:;;;,:::;,:::;;)\qquad
  \gyoung(;;,:;;,:::;,:::;;)\qquad
  \gyoung(;;,:;;,::;;,::;;)
\]
Though a similar definition appears in a lot of papers, nearly all of them do
not contain any other occurrence of the word ``connected''. The definition
that is really used is that consecutive lines overlap exactly at one box.

The second example is the notion of standardization of a word. When you sort a
word $w$, you have to compare the letters, some of which might be equal. If
you want to do a stable sort, you have to consider that left occurrences are
smaller than the right ones. There is exactly one permutation that has the
same relative order of its values. It is usually defined (again from one the
author's paper) as follows:
\begin{quotation}
  With each word $w$ of $A^*$ of length $n$, we associate a permutation
  $\std(w)\in\SG_n$ called the \emph{standardization} of $w$ defined as the
  permutation obtained by iteratively scanning $w$ from left to right, and
  labelling $1,2,\ldots$ the occurrences of its smallest letter, then
  numbering the occurrences of the next one, and so on.
\end{quotation}
Such a process is both very hard to reason with and is very unsuitable for
being implemented in a functional programming language. The good definition is
the following: $\std(w)$ is the unique permutation that has the same
inversions (or non-inversion called here version) as $w$:
\begin{lstlisting}
Definition versions (w : seq T) : rel nat := fun i j =>
  (i <= j < size w) && (nth inh w i <= nth inh w j).
Definition eq_inv (w1 : seq T1) (w2 : seq T2) :=
  (versions w1) =2 (versions w2).
\end{lstlisting}
Then the specification of the standardization map is stated as the following
boolean reflection lemma:
\begin{lstlisting}
  reflect (eq_inv u v) (std u == std v).
\end{lstlisting}
It took nearly three weeks (out of six months to get the first \LR proof)
to clean up this particular very basic notion ending with a thousand lines
long file.

\medskip Though a large part of the library is reusable, it is still very
focused on partitions, tableaux and symmetric functions. We would like to
conclude by assessing the difficulty of writing a widely applicable algebraic
combinatorial library. For this purpose we compare our development with the
files in the \Sage library concerning partitions, permutations, tableaux, skew
tableau, symmetric functions and symmetric group representations. We feel that
this comparison is meaningful since both development are based on a similar
ground: on the combinatorial side the implementation of combinatorial data
structure such as list, partial orders and graphs and on the algebraic side,
polynomials, groups, matrices and representations. We estimate that our 35000
lines of codes formalize less that 5 percent of their counterparts in \Sage. A
rough evaluation shows that the files corresponding to these notions in \Sage
contain around 2300 functions (occurrences of the keyword \verb|def|) and 49k
lines of code.

So having their certified equivalent should need something on the order of a
million lines of code, and the whole combinat library of sage should need 10
million lines. Clearly, a development of this size can only be done if the
community decides to massively invest time and money.

\appendix
\section{Table of content of the library}
Here is a extraction of the table of content of the library thanks to the
\verb|coqdoc| tool. We only gives the details for the most relevant files,
often abreviated and some of the minor files are omitted.
\label{appendix:toc}
\begin{figure*}[h]
  \centering
  \includegraphics[
                   width=\textwidth]{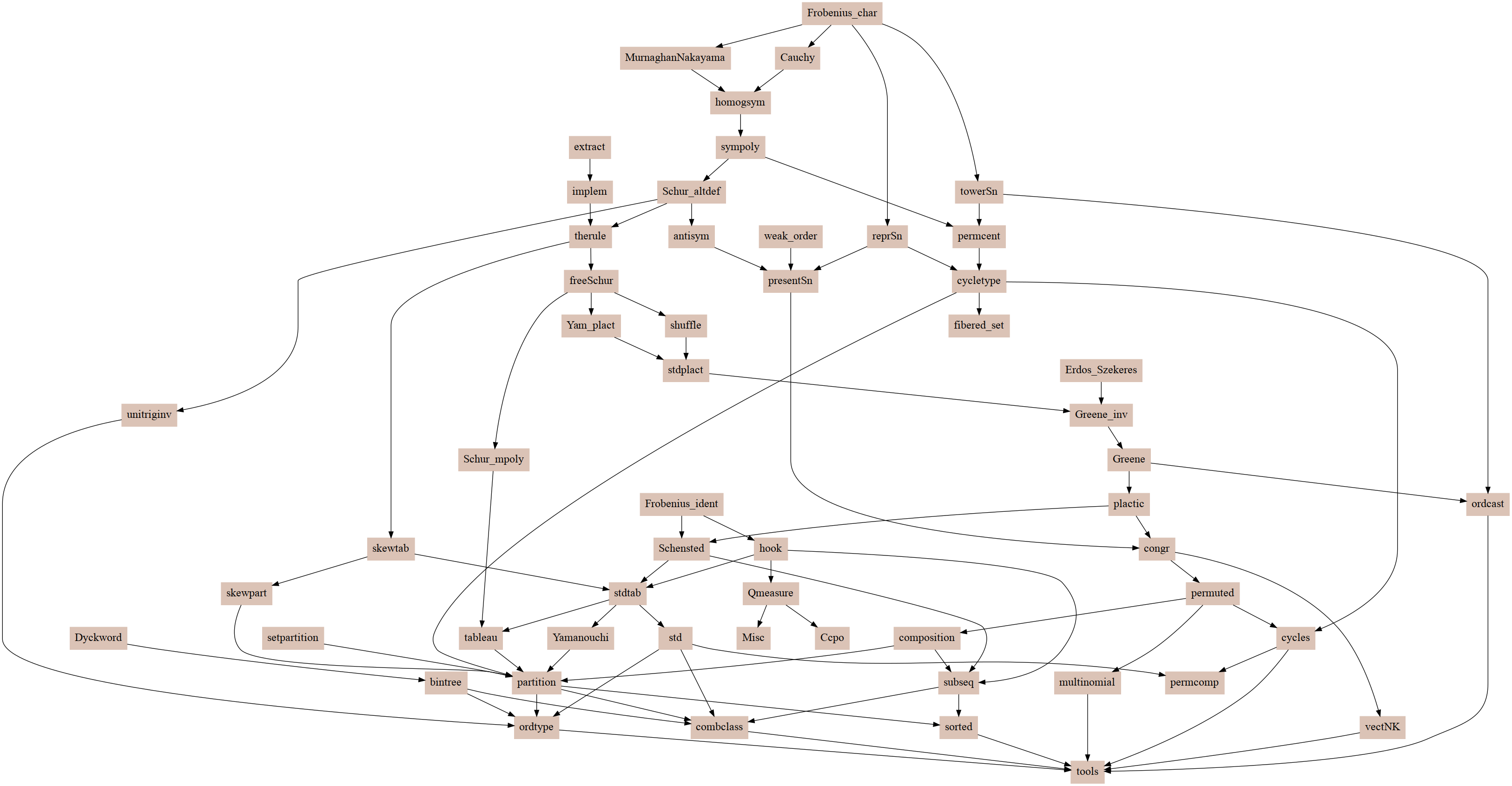}
  \caption{Dependency graph of the files in the library}
  \label{fig.depend}
\end{figure*}

\subsubsection*{Combi.partition: Integer Partitions}
\begin{itemize}
\item Shapes and Integer Partitions
\item Inclusion of Partitions and Skew Partitions
\item Sigma Types for Partitions
\item Counting functions
\item A finite type finType for coordinate of boxes inside a shape
\item The union of two integer partitions
\item Young lattice on partition
\item Dominance order on partition
\item Shape of set partitions and integer partitions
\end{itemize}
\subsubsection*{Combi.setpartition: Set Partitions}
\subsubsection*{Combi.skewpart: Skew Partitions}
\subsubsection*{Combi.skewtab: Skew Tableaux}
\subsubsection*{Combi.std: Standard Words, i.e. Permutation as Words}
\subsubsection*{Combi.stdtab: Standard Tableaux}
\subsubsection*{Combi.subseq: Subsequence of a sequence as a fintype}
\subsubsection*{Combi.tableau: Young Tableaux}
\begin{itemize}
\item Dominance order for rows
\item Tableaux : definition and basic properties
\item Tableaux from their row reading
\item Sigma type for tableaux
\item Tableaux and increasing maps
\end{itemize}
\subsubsection*{Combi.Yamanouchi: Yamanouchi Words}
\subsubsection*{Erdos\_Szekeres.Erdos\_Szekeres: The Erdös-Szekeres theorem}
\subsubsection*{LRrule.extract: Extracting the implementation to OCaml}
\subsubsection*{LRrule.freeSchur: Free Schur functions}
\begin{itemize}
\item Free Schur functions
\item Commutative image of an homogeneous langage
\item Free Schur functions : lifting Schur functions in the free algebra
\item The free Littlewood-Richardson rule
\item Invariance with respect to the choice of the Q-Tableau
\item Conjugating tableaux in the free LR rule
\end{itemize}
\subsubsection*{LRrule.Greene: Greene monotone subsequence numbers}
\begin{itemize}
\item Greene monotone subsequences numbers
\item Greene subsequence numbers
\item Greene number for tableaux
\end{itemize}
\subsubsection*{LRrule.Greene\_inv: Greene subsequence theorem}
\begin{itemize}
\item k-Support and order duality
\item Swaping two letters in a word and its k-supports
\item Greene numbers are invariant by each plactic rules
\item Robinson-Schensted and the plactic monoid
\end{itemize}
\subsubsection*{LRrule.implem: A Coq implementation of the L-R rule}
\subsubsection*{LRrule.plactic: The plactic monoid}
\begin{itemize}
\item Definition of the Plactic monoid
\item Plactic monoid and Robinson-Schensted map
\item Removing the last big letter and plactic congruence
\item Plactic congruence and increasing maps
\end{itemize}
\subsubsection*{LRrule.Schensted: The Robinson-Schensted correspondence}
\begin{itemize}
\item Schensted's algorithm
\item Robinson-Schensted correspondence
\end{itemize}
\subsubsection*{LRrule.shuffle: Shuffle and shifted shuffle}
\subsubsection*{LRrule.stdplact: Plactic congruences and standardization}
\subsubsection*{LRrule.therule: The Littlewood-Richardson rule}
\begin{itemize}
\item Gluing a standard tableaux with a skew tableau
\item Littlewood-Richardson Yamanouchi tableaux
\item The final bijection
\item A slow way to compute LR coefficients in coq
\item The statement of the Littlewood-Richardson rule
\item Pieri's rules
\end{itemize}
\subsubsection*{LRrule.Yam\_plact: Plactic classes and Yamanouchi words}
\subsubsection*{HookFormula.hook: A proof of the Hook-Lenght formula}
\subsubsection*{Basic.congr: Rewriting rule and congruencies of words}
\subsubsection*{SymGroup.cycles: The Cycle Decomposition of a Permutation}
\subsubsection*{SymGroup.cycletype: The Cycle Type of a Permutation}
\begin{itemize}
\item Cycle type and conjugacy classes in the symmetric groups
\item Cycle indicator
\item Central function for 'S\_n
\end{itemize}
\subsubsection*{SymGroup.Frobenius\_char: Frobenius characteristic}
\begin{itemize}
\item Frobenius / Schur character theory for the symmetric groups.
\item Definition and basic properties
\item Frobenius Characteristic and omega involution
\item The Frobenius Characteristic is an isometry
\item The Frobenius Characteristic is a graded ring morphism
\item Combinatorics of characters of the symmetric groups
\item Young characters and Young Rule
\item Irreducible character
\item Frobenius character formula for 'SG\_n
\item Littlewood-Richardson rule for irreducible characters
\end{itemize}
\subsubsection*{SymGroup.permcent: The Centralizer of a Permutation}
\subsubsection*{SymGroup.presentSn: The Coxeter Presentation of the Symmetric Group}
\subsubsection*{SymGroup.towerSn: The Tower of the Symmetric Groups}
\begin{itemize}
\item The Tower of the Symmetric Groups
\item External product of class functions
\item Injection morphism of the tower of the symmetric groups
\item Restriction formula
\item Induction formula
\end{itemize}
\subsubsection*{MPoly.antisym: Antisymmetric polynomials and Vandermonde product}
\begin{itemize}
\item Symmetric and Antisymmetric polynomials
\item Alternating polynomials
\item Vandermonde matrix and determinant
\end{itemize}
\subsubsection*{MPoly.Cauchy: Cauchy formula for symmetric polynomials}
\subsubsection*{MPoly.homogsym: Homogenous Symmetric Polynomials}
\subsubsection*{MPoly.sympoly: Symmetric Polynomials}
\begin{itemize}
\item The Ring of Symmetric Polynomials
\item Multiplicative bases.
\item Littlewood-Richardson and Pieri rules
\item Change of scalars
\item Bases change formulas
\item Fundamental theorem of symmetric polynomials
\item Change of the number of variables
\end{itemize}
\subsubsection*{MPoly.Schur\_mpoly: Schur symmetric polynomials}
\subsubsection*{MPoly.Schur\_altdef: Alternants definition of Schur polynomials}
\begin{itemize}
\item Definition of Schur polynomials as quotient of alternant and Kostka numbers
\item Piery's rule for alternating polynomials
\item Jacobi's definition of schur function
\item Schur polynomials are symmetric at last
\item Kostka numbers
\item Straightening of alternant polynomials
\end{itemize}
\subsubsection*{MPoly.MurnaghanNakayama: Murnaghan-Nakayama rule}

\begin{acks}
  The library contains contributions from Thibaut Benjamin (cycle type of
  permutation and induction of characters of the symmetric groups),
  Jean-Christophe Filliâtre (Why3 implementation of the LR rule), Christine
  Paulin (\MC binding for ALEA + hook length formula), Olivier Stietel (hook
  length formula) and Cyril Cohen (MathComp compatibility + nix CI). The
  authors are very grateful to Georges Gonthier, Cyril Cohen, Pierre Yves
  Strub and Kazuhiko Sakaguchi for helpful discussions. We would like to
  particularly thank Assia Mahoubi and Guillaume Melquiond for their help when
  writing the present paper.

  This project has received funding from the European Research Council (ERC)
  under the European Union's Horizon 2020 research and innovation programme
  (grant agreement No. 101001995).
\end{acks}

\bibliographystyle{ACM-Reference-Format}
\bibliography{lrproof}
\label{sec:biblio}

\end{document}